\documentclass[amsmath,article,aps,fleqn]{revtex4}

\begin{document}

\title{Generalizing Ramanujan's $J$ Functions}
\author{Jerome Malenfant}
\affiliation{American Physical Society\\ Ridge, NY}
\date{\today}
\begin{abstract} We generalize Ramanujan's expansions of the fractional-power Euler functions\\ $  (q^{1/5})_{\infty}  = [~J_1-q^{1/5}+q^{2/5}J_2~](q^5)_{\infty}$ and $  (q^{1/7})_{\infty}  = [~J_1+q^{1/7}J_2 -q^{2/7} +q^{5/7}J_3~] (q^7)_{\infty}$ to $ (q^{1/N})_{\infty} $, where $N$ is a prime number greater than 3.  We show that there are exactly $(N+1)/2$ non-zero $J$ functions in the expansion of $ (q^{1/N})_{\infty} $, that one of these functions has the form  $ \pm q^{X_0}$, that all others have the form $ \pm q^{X_k} \times $ the ratio of two Ramanujan theta functions, and that the product of all  the non-zero $J$'s is $ \pm q^Z$, where $Z$ and the $X$'s denote non-negative integers.
\end{abstract}

\maketitle
\newtheorem{A} {Theorem}
\newtheorem{B}[A]{Theorem}
\newtheorem{I}{ }

\section{Introduction}
In his study of the congruence-5 properties of the partition function $p(n)$,  Ramanujan \cite{Ram} made the replacement
$q \rightarrow q^{1/5}$ in its generating-function equation,
           \begin{eqnarray}
           \sum_{n=0}^{\infty} p(n) q^n = \frac{1}{(q)_{\infty}},
           \end{eqnarray}
where 
            \begin{eqnarray}
           (q)_{\infty} \equiv \prod_{k=1}^{\infty} (1-q^k)
           \end{eqnarray}
is the Euler function.  Then, using Euler's pentagonal number theorem,
           \begin{eqnarray}
           (q)_{\infty} = \sum_{m=-\infty}^{\infty} (-1)q^{m(3m-1)/2},
           \end{eqnarray}
 he made the expansion 
             \begin{eqnarray}
              \frac{ (q^{1/5})_{\infty}}  {(q^5)_{\infty}} = J_1-q^{1/5}+q^{2/5}J_2.
              \end{eqnarray}       
In this equation, the $J$ functions denote power series expansions in $q$ with integer exponents and coefficients.  
These functions can be expressed as the ratios \cite{Somos}
    \begin{eqnarray}
    J_1(q) =   \frac{f(-q^2,-q^3)}{f(-q,-q^4)},~~
    J_2(q) = -~ \frac{f(-q,-q^4)}{f(-q^2,-q^3)},  
    \end{eqnarray}
(sequences A003823 and A007325 in OEIS \cite{integer}, respectively), where $f(a,b)$ is the Ramanujan theta function:
         \begin{eqnarray}
         f(a,b) =f(b,a) \equiv \sum_{n= -\infty}^{\infty} a^{n(n+1)/2}b^{n(n-1)/2} .
        \end{eqnarray} 
 Ramanujan then showed that
             \begin{eqnarray}
             \frac{1}{J_1-q^{1/5}+q^{2/5}J_2 } 
               &=&  \frac{J_1^4+3qJ_2 +q^{1/5}(J_1^3+2qJ_2^2) +q^{2/5} (2J_1^2+qJ_2^3)
             +q^{3/5}(3J_1+qJ_2^4) +5q^{4/5}}{J_1^5 -11q + q^2J_2^5} 
              \end{eqnarray}  
by rationalizing the denominator on the left and using the identity  
               \begin{eqnarray}
               J_1J_2 =-1.
               \end{eqnarray}
From eqs. (1), (7), and another identity, 
               \begin{eqnarray}
              J_1^5 -11q + q^2J_2^5 = \frac{(q)^6_{\infty}}{(q^5)^6_{\infty}} ,
              \end{eqnarray}
it follows that
                \begin{eqnarray}
                \sum_{n=0}^{\infty} p(5n+4)q^n = 5~ \frac{(q^5)^5_{\infty}}{(q)^6_{\infty}} 
                \end{eqnarray}
and therefore that $p(5n+4) \equiv 0  {\rm ~mod~}5$.  
 
In like manner, in studying the congruence-7 properties of $p(n)$, Ramanujan wrote down the expansion 
              \begin{eqnarray}
              \frac{ (q^{1/7})_{\infty}} {(q^7)_{\infty}} = J_1+q^{1/7}J_2 -q^{2/7} +q^{5/7}J_3 
             \end{eqnarray} 
and showed that these $N=7~J$  functions satisfy (among others; see Section IV) the identity
              \begin{eqnarray}
               J_1J_2J_3 = -1.
               \end{eqnarray}

  In this article we generalize these expansions  to $(q^{1/N})_{\infty}$, where $N$ will denote a prime number greater than 3, 
  and we will derive explicit formulas for the $J_p$ functions.    In this, we will be using a slightly different, more convenient notation 
  than that used by Ramanujan, in that the subscript for $J_p$ will correspond to its associated fractional exponent.  I.e., our 
  expansion will read
       \begin{eqnarray}
     \frac{ (q^{1/N})_{\infty}}{(q^{N})_{\infty} }  = J_0 + q^{1/N} J_1  + q^{2/N} J_2  + \cdots   + q^{(N-1)/N} J_{N-1}.
       \end{eqnarray}  
  This equation is equivalent to ``multisecting'' the power series for $(q)_{\infty}$; see  the article by Somos \cite{Somos}.

  \section{Expansion of $(q^N)_{\infty}/(q^{1/N})_{\infty}$}

We make the replacement $q \rightarrow q^{1/N}$ in the Euler function and  write  the identity                                                 
          \begin{eqnarray}
          \frac{1}{(q^{1/N})_{\infty} }   = \frac{1}{(q^{1/N})_{\infty} }
                   \times \frac{ \prod_{p=1}^{N-1} (\omega^p q^{1/N})_{\infty} }{ \prod_{p=1}^{N-1} (\omega^p q^{1/N})_{\infty} }
           =  \frac{ \prod_{p=1}^{N-1} (\omega^p q^{1/N})_{\infty} }{ \prod_{p=0}^{N-1} (\omega^p q^{1/N})_{\infty} }  , 
           \end{eqnarray}
where $\omega \equiv e^{2 \pi i /N} $ is an $N$-th root of unity.   We  consider the product in the denominator in this expression,
with the replacement $x = q^{1/N}$:
                  \begin{eqnarray}                 
                 \prod_{p=0}^{N-1} (\omega^p x)_{\infty} 
                  =  \prod_{p=0}^{N-1} \prod_{k=1}^{\infty}  \left(1- (\omega^{p} x)^k \right). 
                  \end{eqnarray}
Now make the change  of index $k=nN+a$, $1 \leq a \leq N$, to get
                  \begin{eqnarray}                 
                  \prod_{p=0}^{N-1} (\omega^p x)_{\infty}  &=&\prod_{p=0}^{N-1}  \prod_{n=0}^{\infty} 
                  \prod_{a=1}^{N} \left(1- \omega^{pa} x^{nN+a} \right) \nonumber \\
                 &=&\prod_{n=0}^{\infty} \prod_{p=0}^{N-1} \left(1-  x^{nN+N}\right)
                  \prod_{a=1}^{N-1} \left(1- \omega^{pa} x^{nN+a} \right)   \nonumber \\                 
                    &=&\prod_{n=0}^{\infty} \left(1-  x^{nN+N} \right)^N  \prod_{p=0}^{N-1} \prod_{a=1}^{N-1}
                  \left(1- \omega^{pa} x^{nN+a} \right) \nonumber  \\
                   &=& \left(x^{N} \right)^N_{\infty} \prod_{n=0}^{\infty}  \prod_{a=1}^{N-1} \prod_{p=0}^{N-1}
                  \left(1- \omega^{pa} x^{nN+a} \right) .
                  \end{eqnarray}
Since $N$ is prime, $1,\omega^a, \omega^{2a},  \ldots , \omega^{(N-1)a}$ are, for fixed $a$, all distinct, and  so $\{1,\omega^a, \omega^{2a}, \ldots , \omega^{(N-1)a}  \}  = \{ 1, \omega, \omega^2 , \ldots , \omega^{N-1} \} $.  We can therefore  make the replacement $ \omega^{pa} \rightarrow \omega^p$ in the product over $p$, since this amounts to simply a re-ordering of the factors:
                 \begin{eqnarray}
                 \prod_{p=0}^{N-1} (\omega^p x)_{\infty}    &=& \left(x^{N} \right)^N_{\infty} \prod_{n=0}^{\infty}  
                 \prod_{a=1}^{N-1} \prod_{p=0}^{N-1}
                  \left(1- \omega^{p} x^{nN+a} \right) \nonumber \\
                  &=& \left(x^{N} \right)^N_{\infty} \prod_{n=0}^{\infty} \prod_{a=1}^{N-1}   (1- x^{N(nN+a)} )\nonumber  \\
                   &=& \left(x^{N} \right)^N_{\infty} \prod_{n=0}^{\infty} \frac{\prod_{a=1}^{N}   \left(1- x^{N(nN+a)} \right)}
                   {1- (x^{N^2})^{n+1} } 
                     =     \frac{   \left(x^{N} \right)^{N+1}_{\infty}  }{\left(x^{N^2} \right)_{\infty} } ,
                  \end{eqnarray}
 where we've used                 
                  \begin{eqnarray}
                  \prod_{p=0}^{N-1} (1- \omega^{p} X ) =1-X^N  .
                  \end{eqnarray}                  
 And so we have
                  \begin{eqnarray}
                   \prod_{p=0}^{N-1} (\omega^p q^{1/N})_{\infty} =  \frac{   \left(q \right)^{N+1}_{\infty}  }{\left(q^{N}  \right)_{\infty} } 
                  \end{eqnarray}
as the denominator on the right side of eq. (14).  
    
We next consider the numerator in this equation.   We make the replacement $ q^{1/N} \rightarrow \omega^p q^{1/N} $ in (13) and 
take the product of  $(\omega^p q^{1/N})_{\infty}$  over $p=1, \ldots, N-1$:
       \begin{eqnarray}
       \prod_{p=1}^{N-1} (\omega^p q^{1/N})_{\infty} =    (q^{N})_{\infty}^{N-1} \prod_{p=1}^{N-1} \left( J_0 +\omega^p q^{1/N} J_1 
        +\omega^{2p} q^{2/N} J_2 + \cdots + \omega^{(N-1)p} q^{(N-1)/N} J_{N-1} \right).
             \end{eqnarray}
We can use the fact that the product $ \prod_{p=0}^{N-1} \left( x_0 + \omega^p x_1  +\cdots + \omega^{(N-1)p}x_{N-1} \right) $ is equal to the determinant of a circulant matrix:    
     \begin{eqnarray}
    \prod_{p=0}^{N-1} \left( x_0 + \omega^p x_1  +\cdots + \omega^{(N-1)p}x_{N-1} \right) &=&   \left|    \begin{array} {ccccc}
                                                                        x_0 & x_{N-1}  & \cdots & x_2 &x_1 \\ 
                                                                        x_1 & x_0  &\cdots & x_3 &x_2\\
                                                                        x_2 & x_1 & \cdots & x_4& x_3 \\
                                                                       \vdots & \vdots & ~ & \vdots & \vdots \\
                                                                        x_{N-2} & x_{N-3}  &\cdots &x_0 & x_{N-1} \\
                                                                        x_{N-1} & x_{N-2}  &\cdots &x_1 & x_0  \end{array} \right| .
       \end{eqnarray} 
We add columns 1 through $(N-1)$ to the $N$-th columnn and write the determinant as 
       \begin{eqnarray}
       \left|    \begin{array} {ccccc}
         x_0 & x_{N-1}  & \cdots & x_2 &x_1 \\ 
         x_1 & x_0  &\cdots & x_3 &x_2\\
         x_2 & x_1 & \cdots & x_4& x_3 \\
         \vdots & \vdots & ~ & \vdots & \vdots \\
         x_{N-2} & x_{N-3}  &\cdots &x_0 & x_{N-1} \\
         x_{N-1} & x_{N-2}  &\cdots &x_1 & x_0  \end{array} \right|        
        &=&  
         \left|    \begin{array} {ccccc}
         x_0 & x_{N-1}  & \cdots & x_2 & ~x_0+\cdots +x_{N-1} \\ 
         x_1 & x_0  &\cdots & x_3 & ~x_0+\cdots +x_{N-1} \\
         x_2 & x_1 & \cdots & x_4& ~ x_0+\cdots +x_{N-1}  \\
         \vdots & \vdots & ~ & \vdots & \vdots \\
         x_{N-2} & x_{N-3}  &\cdots &~x_0 & x_0+\cdots +x_{N-1}  \\
         x_{N-1} & x_{N-2}  &\cdots &~x_1 & x_0+\cdots +x_{N-1}   \end{array} \right|  \nonumber \\ 
         &=& ( x_0 + \cdots + x_{N-1} )  \left|    \begin{array} {ccccc}
         x_0 & x_{N-1}  & \cdots & x_2 &1 \\ 
         x_1 & x_0  &\cdots & x_3 &1\\
         x_2 & x_1 & \cdots & x_4& 1 \\
         \vdots & \vdots & ~ & \vdots & \vdots \\
          x_{N-2} & x_{N-3}  &\cdots &x_0 & 1 \\
          x_{N-1} & x_{N-2}  &\cdots &x_1 & 1  \end{array} \right| .                                                                                           
          \end{eqnarray}  
And so the product from $p= 1$ to $N-1$ is
        \begin{eqnarray}
       \prod_{p=1}^{N-1} \left( x_0 + \omega^p x_1  +\cdots + \omega^{(N-1)p}x_{N-1} \right) &=& 
         \left|    \begin{array} {ccccc}
         x_0 & x_{N-1}  & \cdots & x_2 &1 \\ 
                                                                        x_1 & x_0  &\cdots & x_3 &1\\
                                                                        x_2 & x_1 & \cdots & x_4& 1 \\
                                                                       \vdots & \vdots & ~ & \vdots & \vdots \\
                                                                        x_{N-2} & x_{N-3}  &\cdots &x_0 & 1 \\
                                                                        x_{N-1} & x_{N-2}  &\cdots &x_1 & 1  \end{array} \right| .
       \end{eqnarray}           
  Upon the  replacements $  x_k \rightarrow  q^{k/N} J_{k}$ we have        
   \begin{eqnarray}
       \prod_{p=1}^{N-1} (\omega^p q^{1/N})_{\infty}  
       =   \left(q^{N}  \right)_{\infty}^{N-1}     \left|    \begin{array} {ccccc}
                                                                        J_0 &  q^{(N-1)/N}J_{N-1}  & \cdots &  q^{2/N}J_2 &1 \\ 
                                                                        q^{1/N}J_1 & J_0  &\cdots & q^{3/N}J_3 & 1\\
                                                                          q^{2/N}J_2 &  q^{1/N}J_1 & \cdots &q^{4/N}J_4&  1\\
                                                                       \vdots & \vdots & ~ & \vdots & \vdots \\
                                                                        q^{(N-2)/N}J_{N-2}& q^{(N-3)/N}J_{N-3}&\cdots &J_0 & 1 \\
                                                                        q^{(N-1)/N}J_{N-1} & q^{(N-2)/N}J_{N-2} &\cdots & q^{1/N}J_1 & 1  \end{array} \right| 
            \end{eqnarray}
 and, together with (14) and (19),       
   \begin{eqnarray}
    \frac{(q^N)_{\infty}}{(q^{1/N})_{\infty} } &=& \frac{1}{ ~ J_0 + q^{1/N} J_1  + q^{2/N} J_2 
    + \cdots + q^{(N-1)/N} J_{N-1} } \nonumber \\~ \nonumber \\
     &=&   \frac{\left(q^{N}  \right)_{\infty}^{N+1} }  { \left(q \right)^{N+1}_{\infty}  }   \left|    \begin{array} {ccccc}
                                                                        J_0 &  q^{(N-1)/N}J_{N-1}  & \cdots &  q^{2/N}J_2 &1 \\ 
                                                                        q^{1/N}J_1 & J_0  &\cdots & q^{3/N}J_3 & 1\\
                                                                          q^{2/N}J_2 &  q^{1/N}J_1 & \cdots &q^{4/N}J_4&  1\\
                                                                       \vdots & \vdots & ~ & \vdots & \vdots \\
                                                                        q^{(N-2)/N}J_{N-2}& q^{(N-3)/N}J_{N-3}&\cdots &J_0 & 1 \\
                                                                        q^{(N-1)/N}J_{N-1} & q^{(N-2)/N}J_{N-2} &\cdots & q^{1/N}J_1 & 1  \end{array} \right|. 
            \end{eqnarray}  

\section{Main Results}  

  \begin{A} Let $N$ be a prime number greater than 3, let $ A$ be an integer $ \in [0, (N-1)/2] $, and let
   \begin{eqnarray*}
      p \equiv \frac{(N-6A)^2-1}{24} {\rm ~mod~} N.
       \end{eqnarray*}
    Then\\~\\
    {\rm (I)}  the expansion
       \begin{eqnarray*}
       \frac{ (q^{1/N})_{\infty}}{(q^{N})_{\infty}}  = J_0(q) + q^{1/N} J_1(q)  + q^{2/N} J_2(q)  + \cdots + q^{(N-1)/N} J_{N-1}(q)
       \end{eqnarray*} 
       has exactly $(N+1)/2$ non-zero terms;  \\~\\
       {\rm (II)}  for  A=0, 
       \begin{eqnarray*}
        J_{p}(q) &=&  (-1)^{  \lfloor (N+1)/6 \rfloor} q^{\lfloor (N^2-1)/24N  \rfloor } ;
       \end{eqnarray*} 
       {\rm (III)}  for A$>$0,   
       \begin{eqnarray*}
        J_{p}(q) &=& (-1)^{ A+ \lfloor (N+1)/6 \rfloor} q^{\lfloor [(N-6A)^2-1]/24N  \rfloor } 
      ~ \frac{f(-q^{2A},-q^{N-2A})}{f(-q^A, -q^{N-A})} .
       \end{eqnarray*}
        Proof:
         \end{A}   
        Prime numbers greater than 3 can be expressed as $N=|6m-1| $ where $m$ is a positive or a negative integer with absolute
        value $|m| = \lfloor (N+1)/6 \rfloor$.\\~\\
    Proof of (I):   We  expand $ (q^{1/N})_{\infty} $ as      
         \begin{eqnarray}
      (q^{1/N})_{\infty}  = \sum_{n=-\infty }^{\infty} (-1)^n q^{n(3n-1)/2N}  .
             \end{eqnarray}  
Set $n = kN+a$, with $ -~ \infty < k < \infty $ and $a=0, \ldots, N-1$.   Then
        \begin{eqnarray}
      (q^{1/N})_{\infty}  &=&  \sum_{a=0}^{N-1} (-1)^a q^{a(3a-1)/2N} \sum_{k=-\infty }^{\infty} (-1)^k q^{k(3kN-6a+1)/2} 
       \end{eqnarray}  
    We now define an equivalence relation on the integers $a \in [ 0 , N -1]$ such that $a_1 \sim a_2 $ iff 
    \begin{eqnarray} 
    \frac{a_1(3a_1-1)}{2} {\rm~ mod~} N \equiv \frac{ a_2(3a_2-1)}{2}  {\rm ~mod~} N .
    \end{eqnarray}  
    We will denote a particular equivalence class either by listing its elements or  as $\{ p \}$, where $p$ is defined by
   \begin{eqnarray} 
   p \equiv   \frac{a(3a-1)}{2} {\rm~ mod~} N 
    \end{eqnarray} 
   for any $a \in \{ p \}$.    From eqs. (13) and (27), each equivalence class corresponds to a non-zero term, with subscript $p$, in the expansion of $  (q^{1/N})_{\infty}$.  I.e.,    
      \begin{eqnarray}
   (q^N)_{\infty} J_p(q) = \left[ (-1)^a  ~ q^{\lfloor a(3a-1)/2N \rfloor } \sum_{k=-\infty }^{\infty} (-1)^k q^{k(3kN-6a+1)/2} \right]_{a \in \{ p\} }   
       \end{eqnarray} 
 where the expression inside the brackets is to be evaluated over the element(s) of the equivalence class $\{ p\}$.
      
   Let $a_1, a_2 \in [0,N-1]$.  From eq. (28),  $a_2 \sim a_1$ iff 
     \begin{eqnarray}
      \frac{(a_2-a_1)[~3(a_2+a_1)-1~]}{2} \equiv 0 {\rm ~mod~}N .
    \end{eqnarray}    
   This requires, since $N$ is prime and $|a_2-a_1|$ is less than $N$, that $N$ divides either $[3(a_2+a_1)-1]$ or $[3(a_2+a_1)-1]/2$, depending on whether  $a_2 + a_1$ is  even or odd.
     \\~\\
    Case 1:  $a_2 + a_1$ is even.  Then 
    \begin{subequations}
    \begin{eqnarray}
    a_2+a_1 = \pm 2Km + \frac{ 1 \mp K}{3} ~~{\rm for~} N = \pm (6m-1)
    \end{eqnarray}
    for some positive integer $K$.  The only solutions for this equation for $a_1,a_2 $ in this interval are:\\~\\
     $a_2  =2m-a_1$ if $m>0$ and $a_1 \leq 2m$;\\
     $a_2 = 2N-|2m| -a_1$ for $m<0$ and $a_1 \geq N-|2m|+1$.\\~\\
     Case 2:  $a_2 + a_1$ is odd.   Then $N$ divides  $(3(a_2+a_1)-1)/2$, and    
      \begin{eqnarray}
    a_2 +a_1 = \pm 4Km + \frac{ 1 \mp 2K}{3}.
    \end{eqnarray}
    \end{subequations}
   The only allowed solution is $a_2 = N+2m-a_1$ for $2m+1 \leq a_1 \leq N+2m$.  
   
     Summarizing, for $m>0$,
        \begin{subequations}   
    \begin{eqnarray}
    a \sim \left\{ \begin{array} {l} 2m-a ~~~~~~~~~~{\rm for~~} a \in [0,2m], \\~\\  N+2m-a ~~~{\rm for~~} a \in [2m+1, N-1] .  \end{array} \right.
    \end{eqnarray}  
    while for $m<0$,  
    \begin{eqnarray}
    a \sim \left\{ \begin{array} {l} N+2m-a ~~~~~{\rm for~~} a \in [0,N-|2m|], \\~\\  2N+2m-a ~~~{\rm for~~} a \in [N-|2m|+1, N-1] .  \end{array} \right.
    \end{eqnarray} 
    \end{subequations}   
     For $m>0$, the first equation is trivial when $a=m$.  Therefore, the equivalence class that contains $m $ has only one distinct element.  
     Similarly, for $m<0$, the class containing $ N +m $ has just one  element.  All other equivalence classes contain exactly 2 elements.  If $M$ is the number of equivalence classes, the $N$ values 
     of $a$ are thus grouped into one 1-element class and $(M-1)$ 2-element classes: $ N =  1 +2(M-1)$.  Therefore, $M=(N+1)/2$, which
     proves (I). \\
       
          Proof of (II):  
    The index $p$ for the 1-element equivalence class, either $\{ a=m\}$ for $N=6m-1$  or $\{ a= N+m\} $ for $N=-6m+1$, is 
     \begin{eqnarray} 
   p \equiv   \frac{m(3am-1)}{2} {\rm~ mod~} N = \frac{N^2-1}{24} {\rm~ mod~} N.
       \end{eqnarray}    
      For $N=6m-1, a=m$, we have from eq. (30) that 
          \begin{eqnarray}
   (q^N)_{\infty}   J_p(q) &=&   (-1)^m q^{\lfloor m(3m-1)/2N \rfloor }    \sum_{k=-\infty }^{\infty} (-1)^k q^{k(3kN-6m+1)/2} \nonumber \\     
          &=&  (-1)^{ \lfloor (N+1)/6 \rfloor} q^{\lfloor (N^2-1)/24N\rfloor }  \sum_{k=-\infty }^{\infty} (-1)^k q^{kN(3k-1)/2}  
            \end{eqnarray}  
 The result then follows, since the sum over $k$ is $(q^N)_{\infty}$ by the pentagonal number theorem. 

The proof for $N=-6m+1,a=N+m$ follows a similar calculation, with the substitution $k =\rightarrow k-1$ in the sum.\\

  Proof of (III):   
    The 2-element equivalence classes are:
     \begin{eqnarray*}
   && m>0:~ \left\{  \begin{array}  {c} ~~~~~~~~~~~~~~~~~~~ 
   \{ 0, 2m \}, ~ \{ 1, 2m -1 \}, ~ \cdots, \{ m-1, m+1\},~~~~~~~~~~~~~~~~~~~~~~~~~~~~~~~~({\rm I})\\  ~\\
    \{ 2m+1, N-1 \}, ~ \{ 2m+2, N-2 \},  ~\cdots,  \{ \frac{1}{2} (N+2m-1),  \frac{1}{2}(N+2m+1) \} ;~~~~~~~~~({\rm II}) \end{array} \right.\\~\\~\\
   && m<0: ~\left\{ \begin{array} {c} ~~~~ \{ 0, N+2m \}, ~ \{ 1, N+2m -1 \},  ~ \cdots,   \{ \frac{1}{2} (N+2m-1), \frac{1}{2} (N+2m+1) \} ,
   ~~~ ~~~~~~~({\rm II }) \\~\\
      \{ N+2m+1, N-1 \}, ~ \{ N+2m+2, N-2 \},  ~\cdots, ~\{ N+m -1, N+m+1\} .~~~~~~~({\rm I})  \end{array}  \right.
     \end{eqnarray*}    
  For a given $m$, they thus break into two groups, which are characterized by the evenness (group I) or the oddness (group II) 
  of  $a_2-a_1$.   To each $\{ a_1,a_2 \},$ with $a_1 < a_2,$ we assign an integer $A$,  defined as 
       \begin{eqnarray}
       A = \left\{ \begin{array} {c}~~~~ (a_2-a_1)/2 ~~~~~~~{\rm if~} a_2-a_1 {\rm ~is~even},  \\~\\  (N-a_2+a_1)/2
       ~~~{\rm if~} a_2-a_1 {\rm ~is~odd}.
        \end{array} \right.
      \end{eqnarray}  
For the 1-element class we set $A=0$. It is easy to see from the list of classes above that each equivalence class corresponds to a different value of $A$ between $0$ and $(N-1)/2$.  
   
  The 2-element classes  thus give 4 cases to consider, which we can characterize as:  \\~\\
    Case 1:  $m>0$, group I:    $~~ a_1 =m-A, ~~~~~~~~~a_2= m+A;$\\   
    Case 2:  $m>0$, group II:   $~ a_1 = m+A,~~~~~~~~~ a_2  = N+m-A;$\\
    Case 3:  $m<0$, group I:    $~~a_1= N+m-A, ~~~a_2 = N+m+A;$\\
    Case 4:  $m<0$, group II:    $ ~a_1 = m+A , ~~~~~~~~~a_2 = N+m-A.$  \\ ~\\    
  Expressed in the variables $N$ and $A$ however,  eq. (29) for $p$ and eq. (30) for $J_p(q)$ takes on the same form in all 4 cases:  
   \begin{eqnarray}
     p &\equiv& \frac{(N-6A)^2-1}{24} {\rm ~mod~} N; \\
    J_p(q)   &=&   (-1)^{A +\lfloor (N+1)/6 \rfloor }  ~ \frac{ q^{\lfloor [(N-6A)^2-1]/24N \rfloor }}{(q^N)_{\infty}}
     \sum_{k=-\infty }^{\infty} (-1)^k q^{Nk(3k-1)/2}  \left[q^{3kA} + q^{(1-3k)A}  \right], ~~(A>0) .
   \end{eqnarray} 
   
Now consider  the identity \cite{theta},
      \begin{eqnarray}
      f(a,b) &=&  (-a;ab)_{\infty} (-b;ab)_{\infty} (ab;ab)_{\infty};~~ (a;q)_{\infty} \equiv \prod _{n=0}^{\infty} (1-aq^n).
      \end{eqnarray}
With this identity, the ratio of the theta functions in part (III) of the theorem is 
       \begin{eqnarray}
      \frac{f(-q^{2A},-q^{N-2A})}{f(-q^A, -q^{N-A})} 
       =  \prod_{n=0} ^{\infty}  \frac{(1-q^{nN+2A})(1-q^{nN+N-2A})}{ (1-q^{nN+A})(1-q^{nN+N-A})  }.
       \end{eqnarray}
We separate the factors in  the numerator into those with even $n$ and with odd $n$:
        \begin{eqnarray} 
        \prod_{n=0} ^{\infty}  \frac{(1-q^{nN+2A})(1-q^{nN+N-2A})}{ (1-q^{nN+A})(1-q^{nN+N-A})  }   
       &=& \prod_{n=0} ^{\infty}  \frac{(1-q^{2nN+2A})(1-q^{(2n+1)N+2A})(1-q^{2nN+N-2A})(1-q^{(2n+1)N+N-2A})}
       { (1-q^{nN+A})(1-q^{nN+N-A})} \nonumber \\
       &=& \prod_{n=0} ^{\infty}  (1+q^{nN+A})(1-q^{2nN+N+2A})(1-q^{2nN+N-2A})(1+q^{nN+N-A}) \nonumber \\    
       &=& \prod_{n=1} ^{\infty}  (1+q^{(n-1)N+A})(1+q^{nN-A}) (1-q^{(2n-1)N-2A})(1-q^{(2n-1)N+2A}) .
       \end{eqnarray}   
 The product on the right, aside from a factor of $(q^N)_{\infty}$, is in the form of the product in the quintuple product
 identity \cite{quintuple}, which we write in the form,
         \begin{eqnarray}
         \prod_{n=1}^{\infty} (1-q^n)(1-aq^{n-1})(1-a^{-1}q^{n})(1-a^2q^{2n-1})(1-a^{-2}q^{2n-1})
         = \sum_{k=-\infty }^{\infty}  q^{k(3k-1)/2} \left[ a^{3k} -a^{1-3k}  \right] ,
       \end{eqnarray}
under the substitutions $q \rightarrow q^N$ and $a \rightarrow -q^A$.  So  we have
     \begin{eqnarray}
      \frac{f(-q^{2A},-q^{N-2A})}{f(-q^A, -q^{N-A})} 
      &=& \frac{1}{(q^N)_{\infty}} \sum_{k=-\infty }^{\infty} (-1)^k q^{Nk(3k-1)/2} \left[ q^{3kA} + q^{(1-3k)A}  \right] .
      \end{eqnarray}
Part (III) of the theorem then follows by comparing the above equation to eq. (38).\\QED.\\

For $N=7$, $A$ takes on the values $0,1,2,3$, corresponding to the functions
\begin{eqnarray*}
J_2  &=& -1,\\
J_0 &=&  ~~\frac{f(-q^2,-q^{5})}{f(-q,-q^6)},\\
J_1 &=&  -  \frac{f(-q^4,-q^{3})}{f(-q^2,-q^5)}, \\
J_5 &=& ~~ \frac{f(-q^{6}, -q,)}{f(-q^3,-q^4)}.
\end{eqnarray*}
The identity (12) then follows trivially.  This identity and that of (8)  can be written (in our notation) as
      \begin{eqnarray}
      J_0J_1J_2 = 1,~~(N=5);~~~
      J_0J_1J_2J_5 = 1,~~(N=7).
      \end{eqnarray}     
 The generalization of these relations is given by the theorem below:     
      
        \begin{A} Let $S = \{p_1, \ldots , p_{(N+1)/2} \}$ be the set of indices 
        corresponding to non-zero $J$ functions in the expansion of $(q^{1/N})_{\infty}$.  Then
         \begin{eqnarray*}
        \prod_{p \in S} J_p(q) = (-1)^{|m|(|m|-1)/2}  q^Z,
         \end{eqnarray*}
where $Z$ is the non-negative integer 
     \begin{eqnarray*}
      Z=\frac{(N-1)(N+1)^2}{48N} - \sum_{p \in S} \frac{p}{N}.
      \end{eqnarray*}
\end{A}
Proof:               
 From parts (II) and (III) of Theorem 1 we have              
        \begin{eqnarray}
       \prod_{ p \in S}  J_p(q) = (-1)^{\sum_A (m+A)}  q^{\sum_A \lfloor [(N-6A)^2-1]/24N \rfloor }   \prod_{ \{a_1,a_2 \} } 
       \frac{f(-q^{2A},-q^{N-2A})}{f(-q^A, -q^{N-A})}  .
               \end{eqnarray}
  The sums over $A$ go from 0 to $(N-1)/2$, while the product on the right is over all 2-element equivalence classes, since the 1-element class contributes only a factor of 1.  Consider the  numerator in this product:
            \begin{eqnarray}
           \prod_{ \{a_1,a_2 \} } f(-q^{2A},-q^{N-2A}) .
           \end{eqnarray}
 There are $(N-1)/2$ factors in this product and each factor contains two distinct positive integers less than $N$; i.e., the exponents $2A$ and $N-2A$, with $A$ ranging from 1 to $(N-1)/2$.  The set of exponents in this product is therefore the set of positive integers less than $N$, and the factors in the product can be reordered as  
        \begin{eqnarray}
        \prod_{ \{a_1,a_2 \} } f(-q^{2A},-q^{N-2A}) =  f(-q,-q^{N-1}) f(-q^{2},-q^{N-2})
        \cdots  f(-q^{(N-1)/2},-q^{(N+1)/2}).
        \end{eqnarray} 
By a similar argument,  the $(N-1)$ exponents in the product in the denominator also equal the set $\{1,2, \ldots , N-1\}$.  The  denominator 
can thus also be reordered as above and cancels with the numerator.\\

The sum over $\lfloor [(N-6A)^2-1] /24N \rfloor$ is found by writing
     \begin{eqnarray}
     \lfloor [(N-6A)^2-1] /24N \rfloor  &=& \frac{ (N-6A)^2-1 }{24N} -\frac{1}{N}~ \frac{ (N-6A)^2-1 }{24}  {\rm~mod~}N
     = \frac{ (N-6A)^2-1 }{24N}  -\frac{p}{N}.
    \end{eqnarray}
We have then
       \begin{eqnarray}
       Z&=& \sum_{A=0}^{(N-1)/2}   \frac{ (N-6A)^2-1 }{24N}  - \sum_{p \in S} \frac{p}{N} 
       = \frac{ (N+1)^2(N-1)}{48N}   -  \sum_{p \in S}   \frac{p}{N} .
       \end{eqnarray}

To find the exponent of (-1), we consider the cases $m>0$ and $m<0$ separately: \\~\\
$m>0$:  $A$ goes from 0 to $(N-1)/2 = 3m-1$;
      \begin{subequations}
     \begin{eqnarray}
    \sum_{A=0}^{3m-1} (m+A)  = m(3m) + \frac{(3m-1)(3m)}{2}  = 3~ \frac{5m^2-m }{2}  .
   \end{eqnarray}  
 But $  (-1)^{3(5m^2-m)/2} = (-1)^{m(m-1)/2}$. \\~\\
$m<0$:  $A$ goes from $0$ to $|3m|$; 
     \begin{eqnarray}
   \sum_{A=0}^{|3m|} (m+A)  = m(|3m|+1) + \frac{|3m|(|3m|+1)}{2}  =  \frac{3m^2-m}{2} .
   \end{eqnarray}  
     \end{subequations}
 In this case,  $ (-1)^{(3m^2-m)/2} = (-1)^{m(m+1)/2}$.  Therefore,  both cases are covered by the factor $ (-1)^{|m|(|m|-1)/2} $.
 
       QED.\\
            
   \section{Some Additional Remarks }   
   
   To derive the identities in eqs. (8) and (9), Ramanujan cubed both sides of eq. (4), used Jacobi's identity,
   \begin{eqnarray}
   (q)^3_{\infty} = \sum_{n=0}^{\infty} (-1)^n (2n+1)q^{n(n+1)/2},
   \end{eqnarray}
   to expand the left-hand side in fractional powers of $q$, and then equated terms.   Another way of arriving at eq. (9) is to use eq. (19) and express the product on the left side as the 
   determinant of a circulant matrix as in eq. (21).  Setting $N=5$ in eq. (19) and dividing by $(q^5)^5_{\infty}$,   we have        
    \begin{eqnarray}
    \frac{   \left(q \right)^{6}_{\infty}  }{\left(q^{5}  \right)^6_{\infty} }  = \frac {1} { (q^5)^5_{\infty} } \prod_{p=0}^{4} (\omega^p q^{1/5})_{\infty}
    &=&\left|    \begin{array} {ccccc} 
          J_0 & 0 &0&q^{2/5}J_2&  -q^{1/5}  \\ 
          -q^{1/5} & J_0 &0 &0 &q^{2/5}J_2\\           
          q^{2/5}J_2 & -q^{1/5} & J_0& 0& 0 \\
         0 & q^{2/5}J_2 & -q^{1/5} & J_0 & 0 \\ 
         0 & 0 &q^{2/5}J_2 &- q^{1/5} & J_0 \end{array} \right|\nonumber \\
         &=&  J_0^5  +q(  5J_0J_2 -1 - 5J_0^2 J_2^2)  + q^2J_2^5 
        \end{eqnarray} 
Now substituting $J_0J_2=-1$ into this equation gives the identity in (9).

Clearly, we can continue in this fashion.  E.g.,  for $N=7$,  this becomes
    \begin{eqnarray}
    \frac{   \left(q \right)^{8}_{\infty}  }{\left(q^{7}  \right)^8_{\infty} }  =  \frac{ 1  }{\left(q^{7}  \right)^7_{\infty} } \prod_{p=0}^{6} (\omega^p q^{1/7})_{\infty}
    &=& \left|    \begin{array} {ccccccc} 
          J_0 & 0 &  q^{5/7}J_5  & 0&0&  -q^{2/7} &  q^{1/7}J_1  \\ 
          q^{1/7}J_1 & J_0 &0 & q^{5/7}J_5&0& 0& -q^{2/7}  \\           
          - q^{2/7} & q^{1/7}J_1 & J_0& 0&  q^{5/7}J_5&0& 0  \\
         0 &- q^{2/7} &q^{1/7}J_1& J_0 & 0 &q^{5/7}J_5 & 0   \\ 
         0 &0& - q^{2/7}& q^{1/7}J_1& J_0 & 0  & q^{5/7}J_5\\
         q^{5/7}J_5  &0&  0 &- q^{2/7}  & q^{1/7}J_1 & J_0 &0 \\
         0 & q^{5/7}J_5  &0& 0  & - q^{2/7} &  q^{1/7}J_1&J_0  \\
         \end{array} \right| \nonumber \\~ \nonumber\\
         &=& J_0^7 +q( ~J_1^7 +7J_0J_1^5 +14J_0^2J_1^3 +7J_0^4J_1^2J_5 +7J_0^3J_1 +7J_0^5J_5  ~  )  \nonumber \\
         \nonumber~\\&&
         +q^2( ~ 7J_0J_1^4J_5^2 +7J_1^3J_5 +7J_0^2J_1^2J_5^2 +14J_0J_1J_5 +14J_0^3J_5^2 -1~)  \nonumber
                  \\ \nonumber~\\ && +q^3(~14J_1^2J_5^3 +7J_0^2J_1J_5^4 +7J_0J_5^3~) +7q^4J_1J_5^5 +q^5J_5^7 .
            \end{eqnarray} 
Using $J_0J_1J_5=-1$, this simplifies to 
    \begin{eqnarray}
    \frac{   \left(q \right)^{8}_{\infty}  }{\left(q^{7}  \right)^8_{\infty} } 
             &=& J_0^7 +qJ_1^7 +q^5J_5^7  +7q( ~ J_0J_1^5  +J_5J_0^5  +q^3J_1J_5^5~ ) 
         +14q( ~J_0^2J_1^3+qJ_5^2 J_0^3 +  q^2J_1^2J_5^3 ~)  -8q^2 .        
            \end{eqnarray}         
We can further simplify this expression by using some of  the other Ramanujan $N=7$ identities \cite{Ram}:
           \begin{subequations}
          \begin{eqnarray}
            \frac{J_0^2}{J_5} + \frac{J_1}{J_5^2}&=& q;  \\
           J_0^7 +qJ_1^7 + q^5J_5^7  &=& \frac{(q)^8_{\infty}}{(q^7)_{\infty}^8} + 14q \frac{(q)^4_{\infty}}{(q^7)_{\infty}^4}
           +57q^2;\\
           J_0^3J_1 +q J_1^3J_5 + q^2J_5^3 J_0  &=& - \frac{(q)^4_{\infty}}{(q^7)_{\infty}^4}-8q ;\\ 
            J_0^2J_1^3 +q J_5^2 J_0^3+q^2 J_1^2J_5^3  &=& - \frac{(q)^4_{\infty}}{(q^7)_{\infty}^4} -5q.  
           \end{eqnarray}
           \end{subequations}
[Note that we have corrected a misprint in \cite{Ram} in the last term on the right in eq. (55b).]  Substituting the left-hand sides of eqs. (55b) and (55d) into eq. (54),  we get the additional identity,
           \begin{eqnarray}
           J_0J_1^5  +J_5J_0^5  +q^3J_1J_5^5  =3q.    
          \end{eqnarray}

\end{document}